\documentclass[graybox]{svmult}

\usepackage{type1cm}        
\usepackage{makeidx}         
\usepackage{graphicx}        
\usepackage{multicol}        
\usepackage[bottom]{footmisc}

\usepackage{newtxtext}       %
\usepackage[varvw]{newtxmath}       

\makeindex             

\usepackage{dsfont}
\usepackage{pgffor}
\usepackage{booktabs}
\usepackage{subcaption}
\usepackage{algpseudocode}

\makeindex             

\begin{document}

\title*{Design of a checkerboard counterflow heat exchanger for industrial applications}
\titlerunning{Design of a checkerboard counterflow heat exchanger } 
\author{Nicola Parolini, Vanessa Covello, Alessandro Della Rocca, and Marco Verani}
\authorrunning{N.~Parolini, V.~Covello, A.~della Rocca, M.~Verani}
\institute{Nicola Parolini \at MOX, Dipartimento di Matematica, Politecnico di Milano, P. Leonardo da Vinci, 32, 20133 Milano \email{nicola.parolini@polimi.it}
\and Vanessa Covello \at MOX, Dipartimento di Matematica, Politecnico di Milano, P. Leonardo da Vinci, 32, 20133 Milano \email{vanessa.covello@polimi.it}
\and Alessandro Della Rocca \at TENOVA, Via Gerenzano, 58, 21053 Castellanza (VA) \email{alessandro.dellarocca@tenova.com}
\and Marco Verani \at MOX, Dipartimento di Matematica, Politecnico di Milano, P. Leonardo da Vinci, 32, 20133 Milano \email{marco.verani@polimi.it}}
%
%
\maketitle

\abstract*{
This work is devoted to the design of a checkerboard air-gas heat exchanger suitable for industrial applications. The design of the heat exchanger is optimized in order to obtain the maximum increase of the outlet air temperature, considering different geometrical design parameters and including manufacturing constraints. The heat exchanger efficiency has been assessed by means of the $\epsilon$-NTU method. The perfomances are compared with traditional finned recuperators and appreciable enhancement of the exchanger efficiency has been observed adopting the new design.} 

\abstract{
This work  is devoted to the design of a checkerboard air-gas heat exchanger suitable for industrial applications. The design of the heat exchanger is optimized in order to obtain the maximum increase of the outlet air temperature, considering different geometrical design parameters and including manufacturing constraints. The heat exchanger efficiency has been assessed by means of the $\epsilon$-NTU method. The perfomances are compared with  traditional  finned recuperators and appreciable enhancement of the exchanger efficiency has been observed adopting the new design.} 

\section{Introduction}
\label{intro}

Heat exchangers are widely used in many industrial area and represent a field of research deeply investigated during the last decades, see, e.g.,  \cite{alam2018, florides2007,srimuang2012}. Their importance has recently gained increasing attention due to the impact on energy conservation, conversion and recovery. An efficient heat exchanger design can affect the entire industrial processes \cite{shah2003}, and a key role is played for heat recovery in high temperature industrial systems, as in burners for industrial furnaces. In that context, a heat exchanger is coupled to the burner and used for pre-heating the combustion air with the exhaust gas, by means of recuperative burners. Several studies have been performed on recuperative burners and their effect on combustion efficiency and NOx production; see, e.g., \cite{wunning1997, galletti2007,demayo2002}. Various types of heat exchangers are currently available and different strategies are used to enhance their performance, see, e.g., the adoption of porous structures, fins or pins of various kinds and shapes, swirled flows, as in \cite{yu1999,kahalerras2008,nagarani2014}. Besides these strategies, the adoption of parametric studies and optimization methods coupled with CFD simulations has emerged as powerful tools to accurately drive the heat exchanger design toward the desired performance, and several examples can be found in the literature \cite{chong2002,rao2010,selbacs2006,selma2014,cavazzuti2015}.  
In this context, the most typical goals include the maximization of the heat transfer and the minimization of the pressure drop, see, e.g., \cite{foli2006, hilbert2006, siavashi2018, mohammadi2020,  amini2014, chamoli2017, abbasi2020,hsieh2012}.                         
 
Starting from the know-how developed in the context of optimization and finite-volume approximation of complex industrial flows \cite{bruggi2011, bruggi2018, NPV2024, negrini2024immersedboundarymethodpolymeric}, in this paper we present an extensive parametric numerical investigation carried out by means of the OpenFOAM open-source library \cite{Weller1998}  for the design of a checkerboard air-gas recuperator based on a recent Tenova patent \cite{dellarocca2018b}, in order to optimize its heat exchange performances. The shape of the exchanger has been properly parametrized and an extensive simulation campaign on various geometrical configurations has been performed exploiting a fully automatic mesh morphing strategy. The heat exchanger efficiency has been assessed by the $\epsilon$-NTU method and its perfomances has been compared with those of existent traditional double pipe finned recuperator. A preliminary assessment based on recent experimental tests is also presented.

The paper is organized as follows. In Section \ref{geom} we present the geometrical configuration of the device and the parametrization that is used to improve its design; in Section \ref{gov-eq}, the governing equations defining the conjugate heat transfer problem are introduced; the computational setup and the results of the simulation campaign are presented in Section \ref{num-res}, while a brief discussion on the experimental validation is given in Section \ref{validation}. Finally, some conclusions are drawn in Section \ref{conclusion}. 

\section{Geometrical parametrization}
\label{geom}

\begin{figure}[h!]
    \centering
     \includegraphics[width=0.6\textwidth]{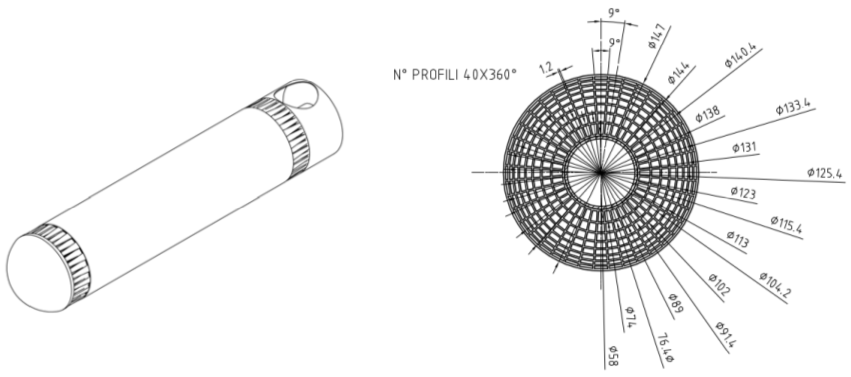}
    \caption{Tenova patent for the recuperative burner \cite{dellarocca2018b}}
    \label{fig:pat}
\end{figure}

The checkerboard arrangement of the heat exchanger under investigation has been developed on the basis of the Tenova patent \cite{dellarocca2018b} displayed in Fig.\ref{fig:pat}.
The parametric analysis that is carried out in this work focuses on the core part of the recuperator, in which heat transfer between air and exhaust gases occurs. It consists of $n \times m$ channels in counterflow mode, arranged in a checkerboard pattern, where $n$ and $m$ denote the radial and azimuthal number of channels, respectively.

\begin{figure}[h]
    \centering
     \includegraphics[width=0.22\textwidth]{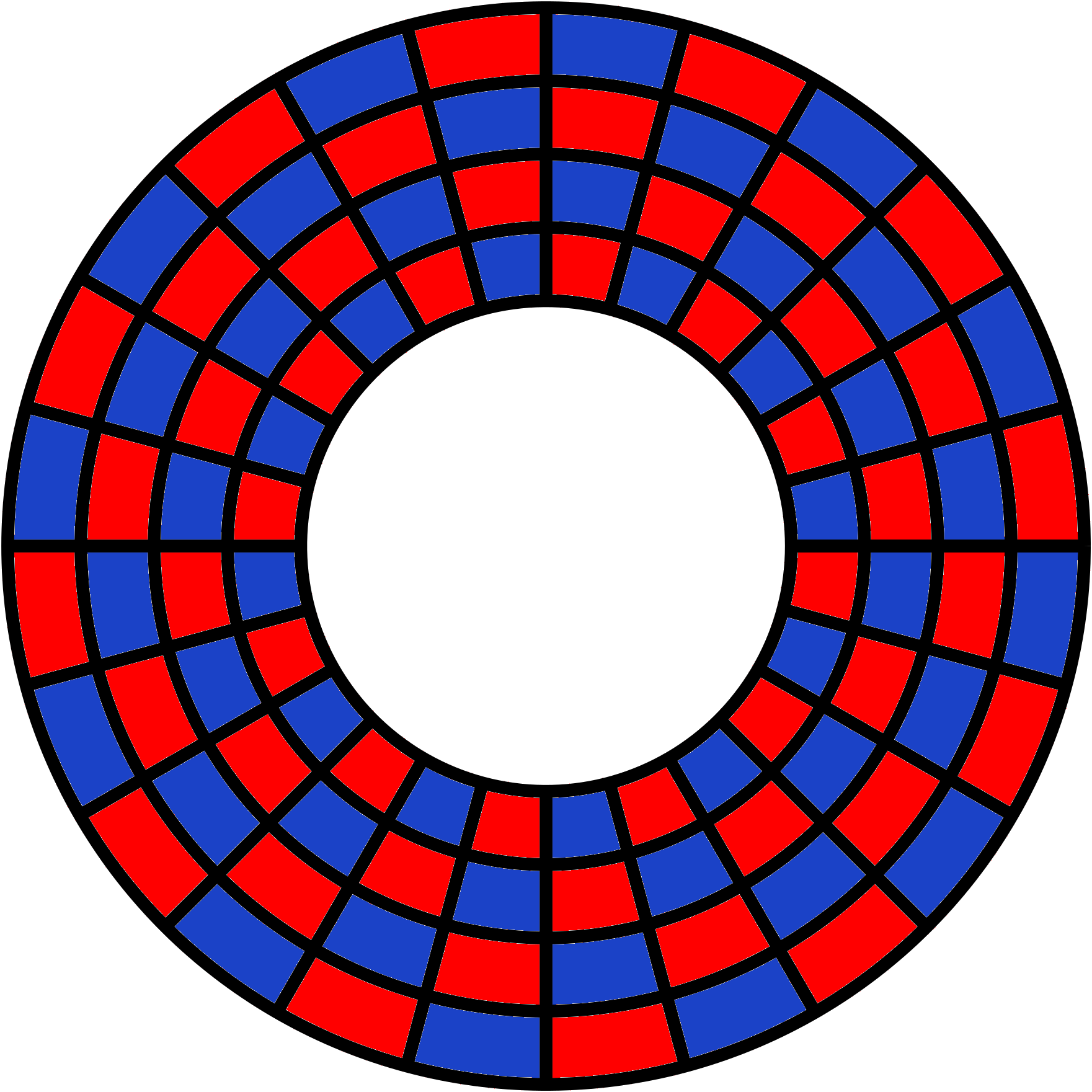}
     \hspace{2mm}
          \includegraphics[width=0.22\textwidth]{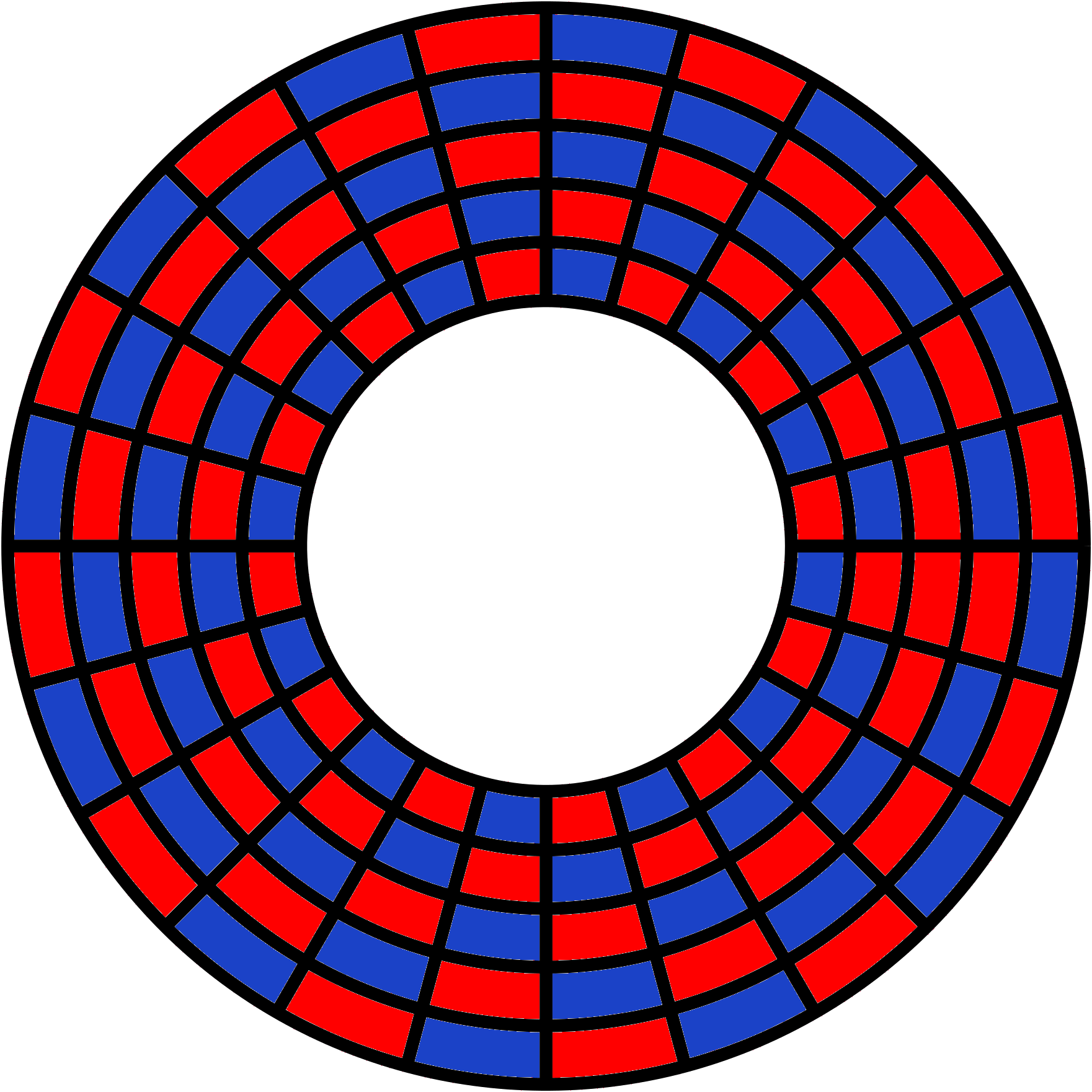}
     \hspace{2mm}
          \includegraphics[width=0.22\textwidth]{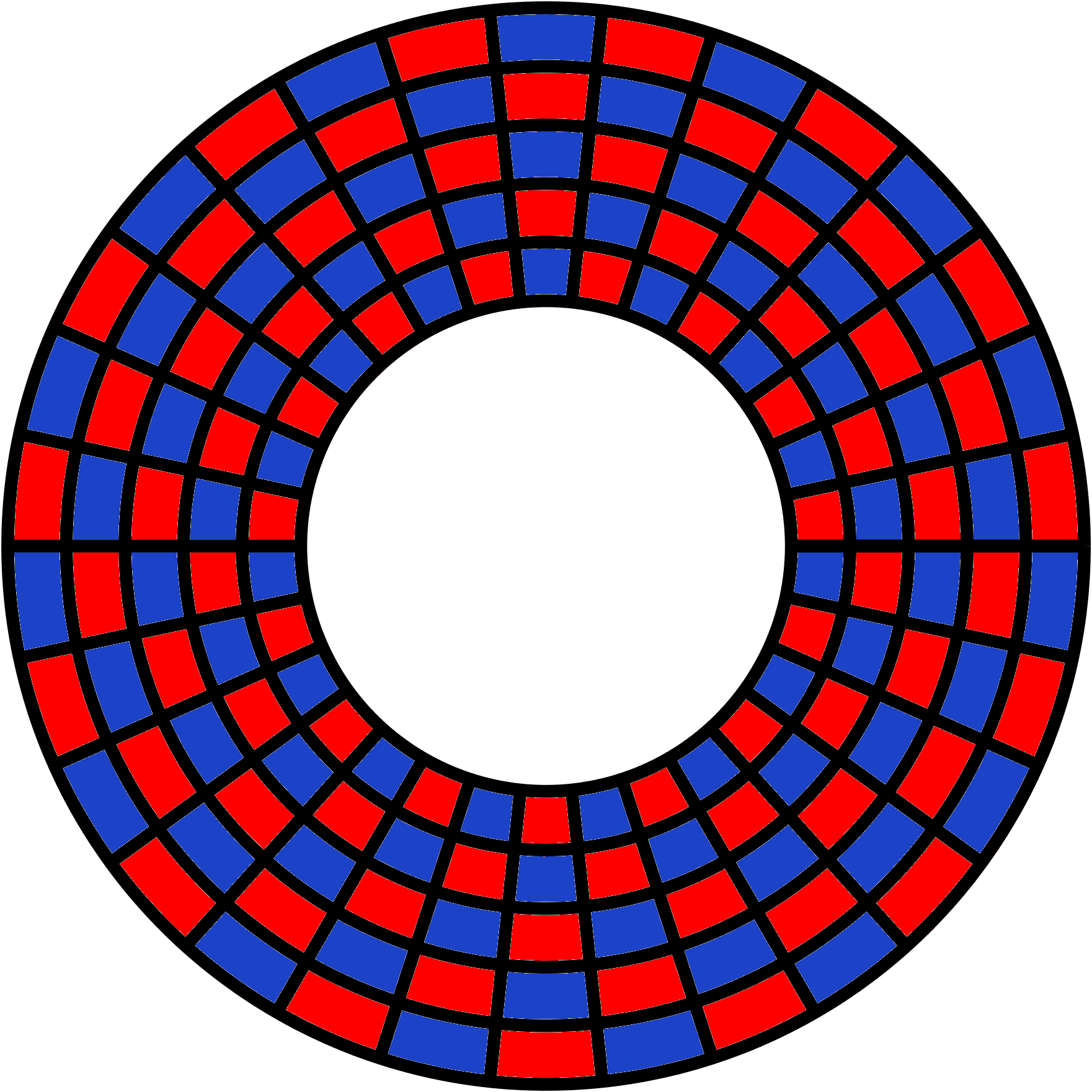}
     \hspace{2mm}
           \includegraphics[width=0.22\textwidth]{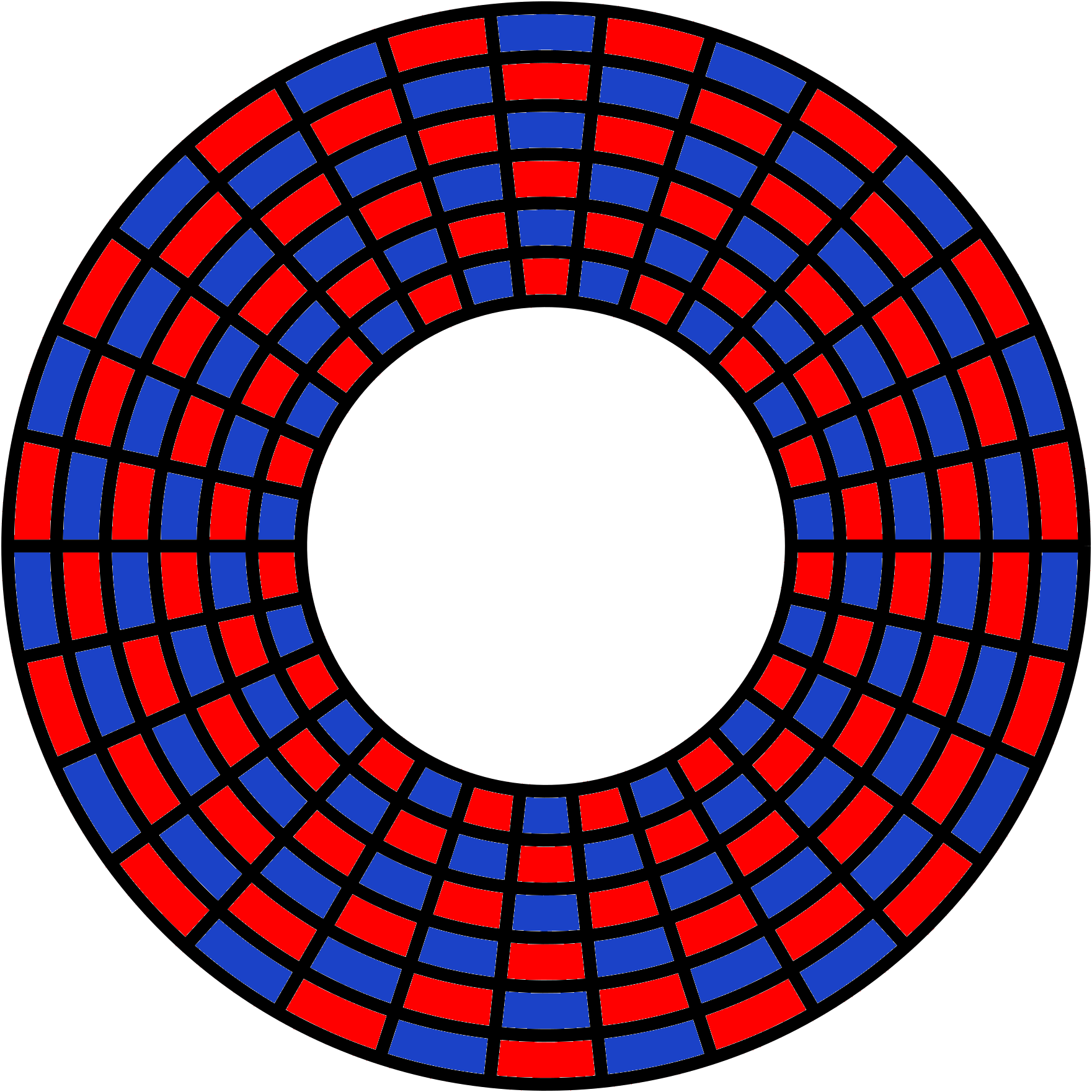}
    \caption{Different configurations of the recuperative burner: from left to right 4x24, 5x24, 5x30, and 6x30 channels. }
    \label{fig:chan}
\end{figure}

\begin{figure}[b!]
    \centering
     \includegraphics[width=0.3\textwidth]{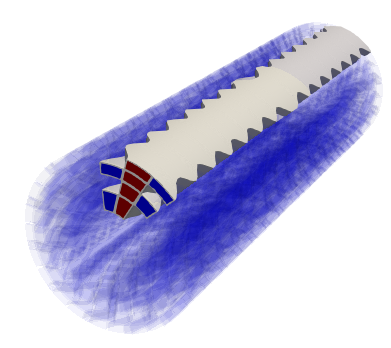}
     \includegraphics[width=0.3\textwidth]{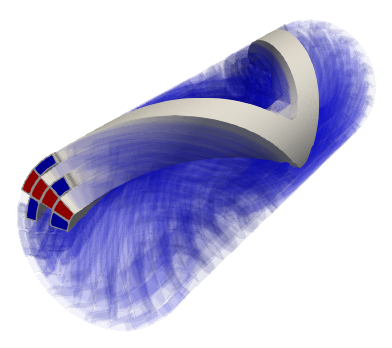}
     \includegraphics[width=0.3\textwidth]{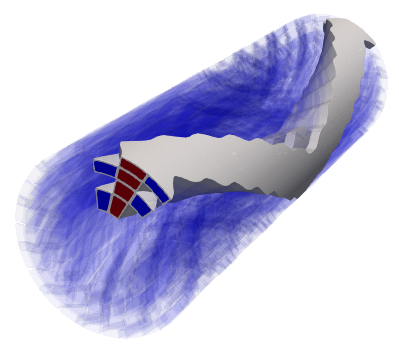}
    \caption{Different geometric perturbations: wavy (top), twist (middle), hybrid twist-wavy (bottom), for the 4x24 channels configurations.}
    \label{fig:geom}
\end{figure}

In particular, we have considered the 4x24, 5x24, 5x30, 6x24 channel arrangements displayed in Figure \ref{fig:chan}. The ranges on the number of radial and azimuthal elements have been selected in order to not exceed the limit of about 5 mm for the minimum channel thickness, with the purpose of guaranteeing the manufacturability and maintainability of the recuperator.  For each channel arrangement, different configurations have been generated by applying a set of geometrical deformations to a baseline configuration with axially straight channels. A sketch of the transformations initially considered, namely a sinusoidal (wavy), a twist and a combined twist-wavy axial transformations, are depicted in Figure \ref{fig:geom}. 

\begin{figure}[t]
    \centering
     \includegraphics[width=0.22\textwidth]{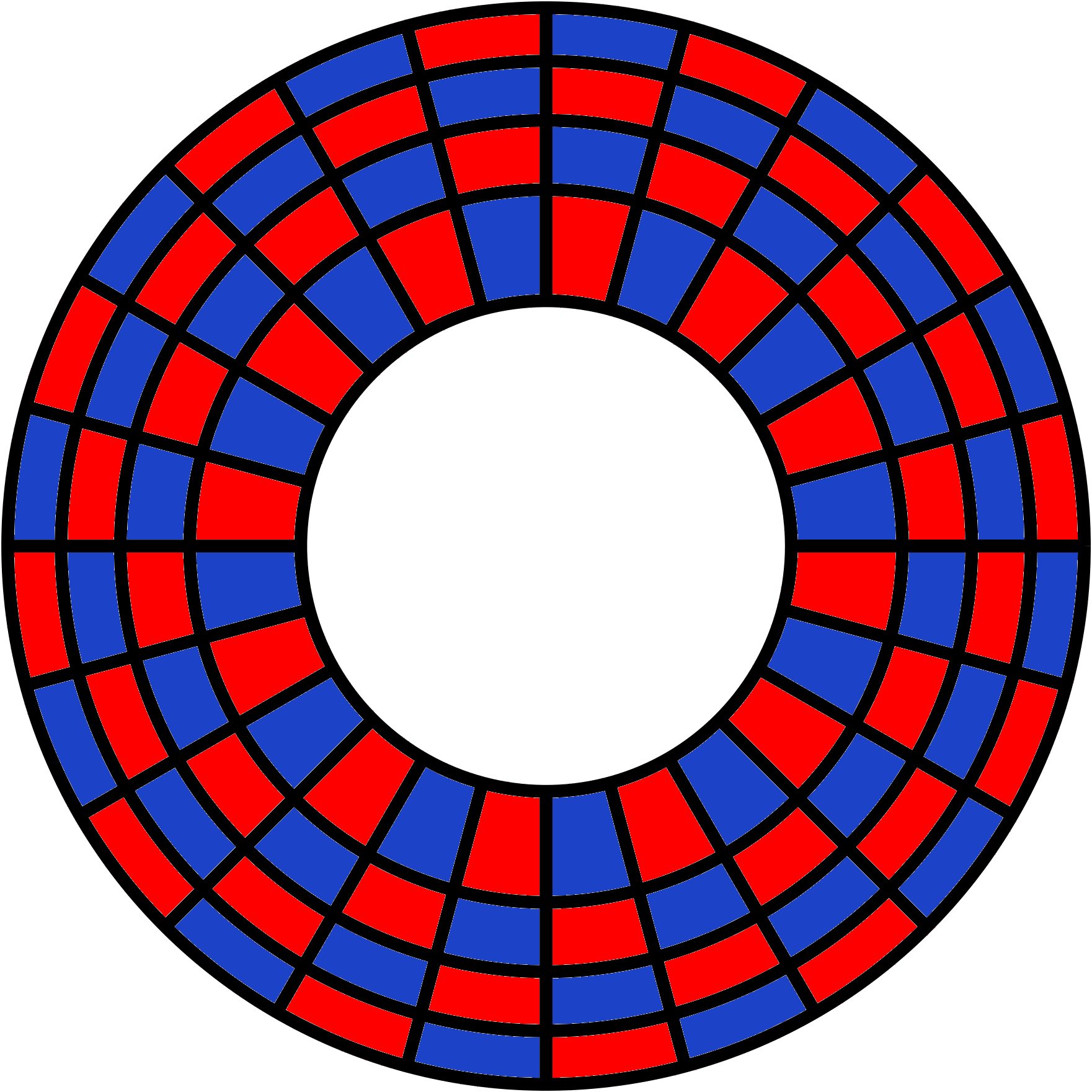}
     \hspace{5mm}
     \includegraphics[width=0.22\textwidth]{img/CheckerSection4x24.png}
     \hspace{5mm}
     \includegraphics[width=0.22\textwidth]{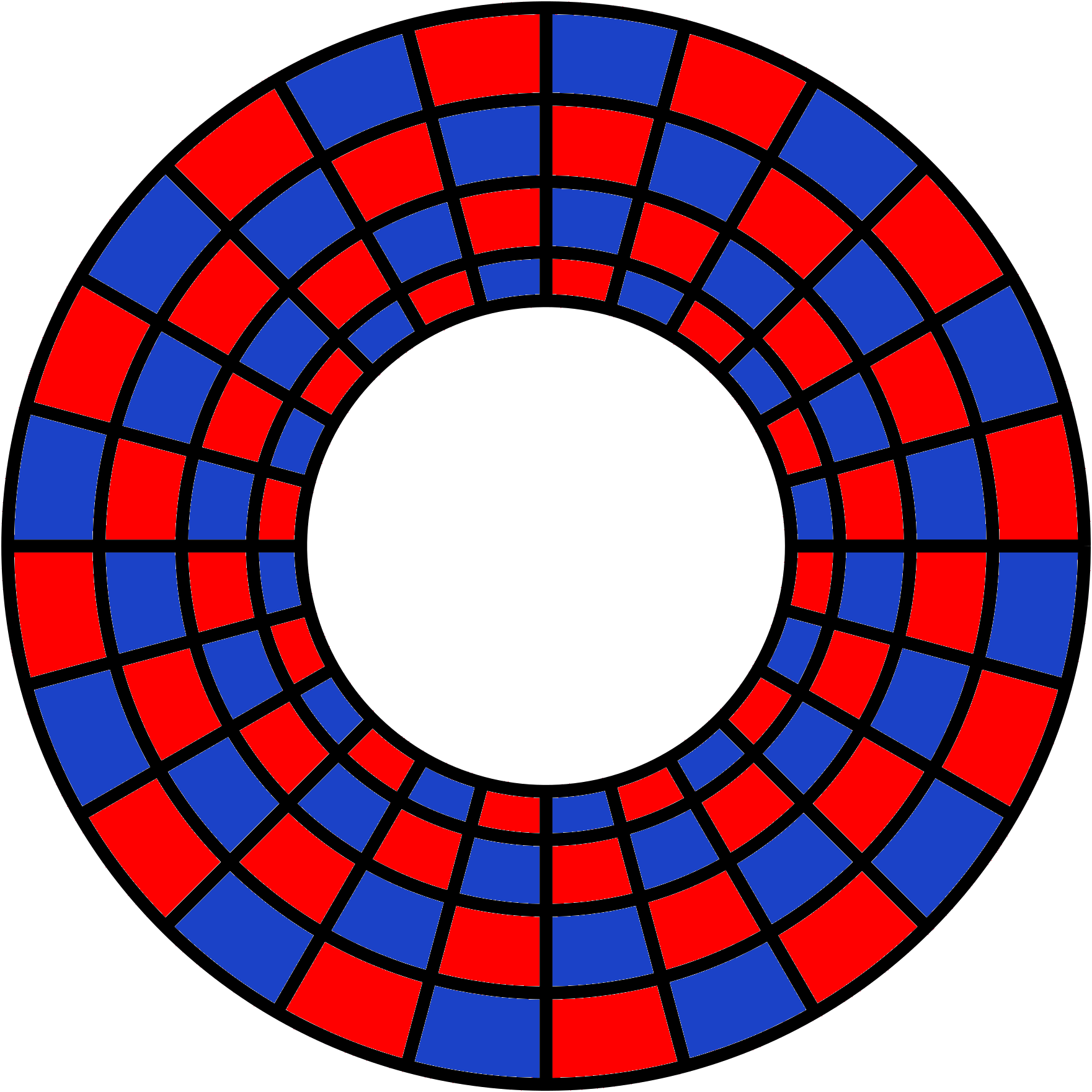}
    \caption{Different scaling perturbation in 4x24 channels configuration: $q<1$ (left), $q=1$ (middle), $q>1$ (right). }
    \label{fig:scal}
\end{figure}

A further scaling geometrical perturbation has been applied for each deformed configuration, in order to scale the channel thickness along the radial direction, as shown in Figure \ref{fig:scal} for the 4x24 case. The radial scaling is driven by a scaling factor $q$, that corresponds to 1 in case of uniform radial thickness, and is lower or higher than 1 if the radial thickness is progressively reduced close to the external or internal diameter, respectively. In our study, a range of admissibility for $q$ has been imposed for each checkerboard layout, in order to guarantee the manufacturability of the recuperator.

\begin{figure}[h]
    \centering
    \includegraphics[width=0.3\textwidth]{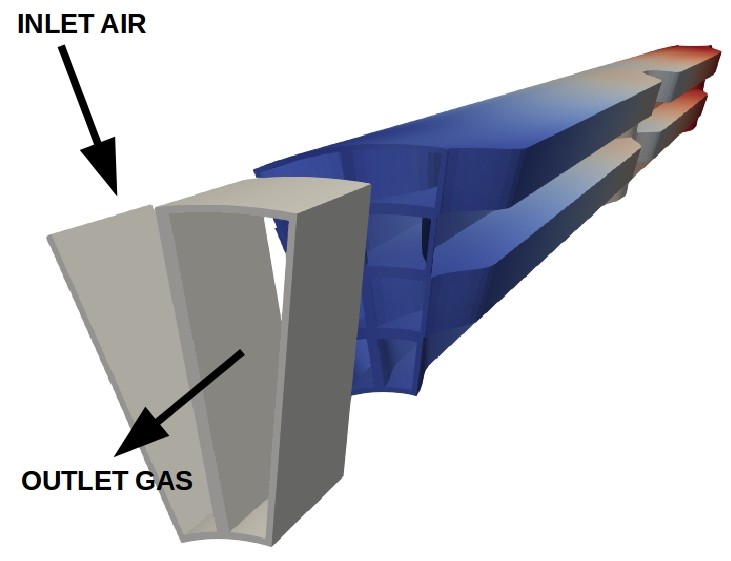}
    \includegraphics[width=0.3\textwidth]{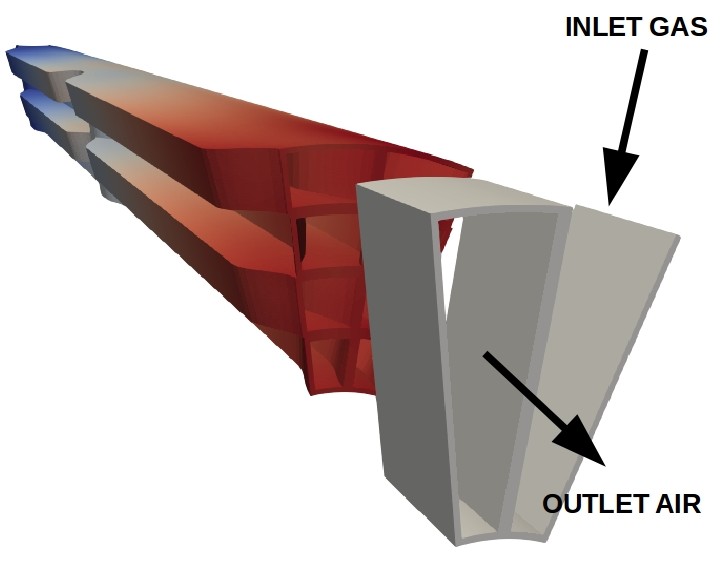}
    \includegraphics[width=0.33\textwidth]{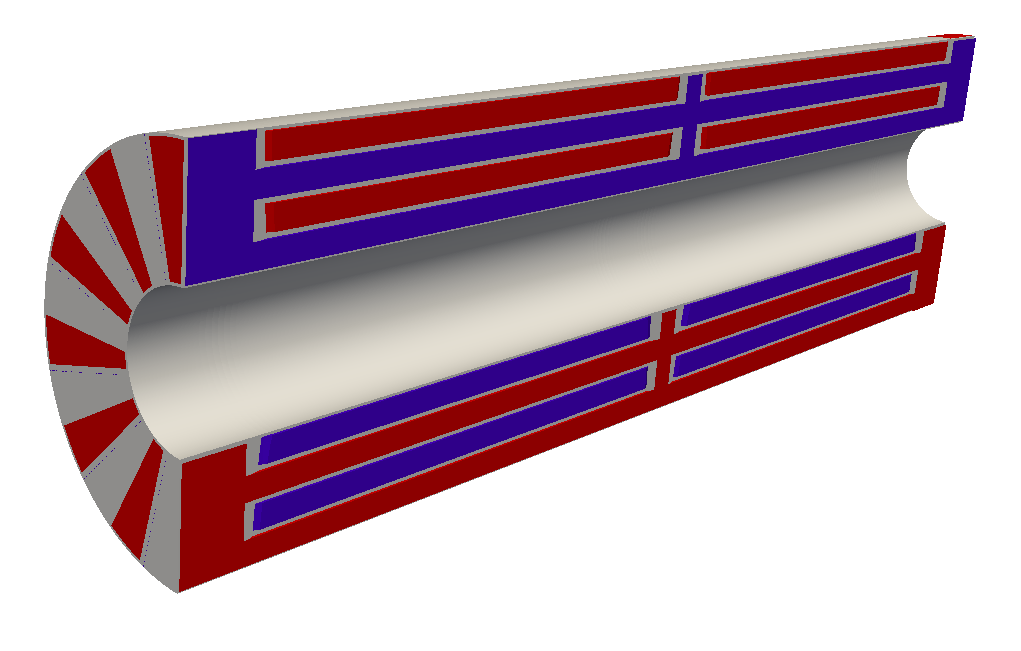}

    \caption{Integration of cold side (left) and hot side (middle) flow collectors and transversal section of the recuperator (right).}
    \label{fig:coll}
\end{figure}

Some constraints have been imposed to the design of the recuparator, related to the additive manufacturing technology that is considered for its production. In particular, the recuperator length  (around 55 cm) exceeds the maximum printing height, therefore the recuperator is printed in two halves.
The 2 checkerboarded cores of the recuperator have been coupled with the air and gas collectors displayed in Figure \ref{fig:coll} (left and middle) collecting the hot and cold channels of each angular sectors in a unique radial inlet and a unique axial outlet.
A similar junction is present at the middle of the recuperator, allowing an easier sealing between the two halves  (see Figure \ref{fig:coll}, right). From now on, the undeformed configuration will be denoted as \textit{Baseline} configuration. 
Additional constraints related to the additive manufacturing process, such as a minimal thickness of 2mm for the solid walls and a maximum overhang angle of 43 degrees has been considered and accounted for by ad-hoc fully automatic mesh generation and mesh morphing strategies.

\begin{figure}[h!]
    \centering
     \includegraphics[width=0.15\textwidth]{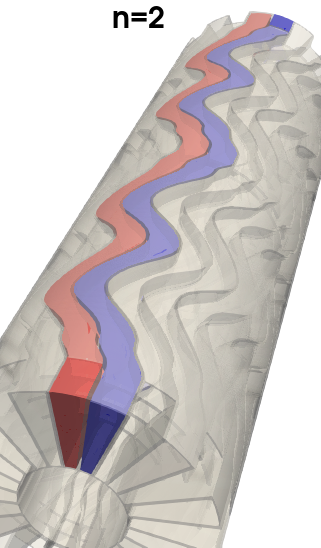}
     \hspace{10mm}
     \includegraphics[width=0.15\textwidth]{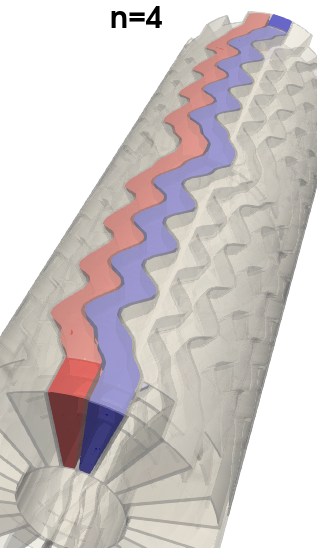}
     \hspace{10mm}
     \includegraphics[width=0.15\textwidth]{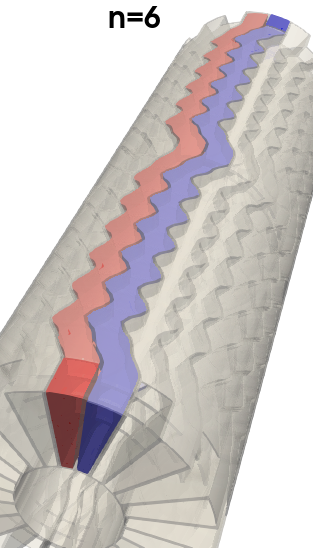}
    \caption{Wavy channels configurations for wavenumber $n$=2 (right), $n$=4 (middle), and $n$=6 (right).}
    \label{fig:wav}
\end{figure}

In this paper, we only report the outcomes of the analyses referred to the wavy and scaling transformations, since the twist and the hybrid twist-wavy transformations were found to be much less effective.
For the wavy transformation, we considered a wave number parameter $n$, ranging from $n$=1 to $n$=6, where $n$ represents the number of sinusoidal oscillations on the half-length of the recuperator, as displayed in Figure \ref{fig:wav}.

The geometric law for the wavy transformation is:
\begin{equation}\label{eq:wavy}
    f(z)=A(\cos\left(2\pi n \frac{z}{L}\right) -1)
\end{equation}
where $A$ denotes the amplitude of the geometrical oscillation and $z$ the axial coordinate and $L$ the axial length of the core portion in each half of the recuperator. Depending on $n$, the amplitude $A$ has been limited to guarantee that the overhang angle, that is the inclination of the print wall from the vertical axis, does not exceed 43 degrees.

\section{Conjugate Heat Transfer model}
\label{gov-eq}
The governing equations are based on the so-called conjugate heat transfer paradigm, widely adopted in the literature \cite{patankar2018, zudin2016, zhang2013, dorfman2010} for this type of analyses. In conjugate heat transfer problems, conservation equations are solved for both fluid and solid regions, and the suitable thermal coupling conditions are imposed at the solid-fluid interface. 
The governing equations expressing mass and momentum balance for the fluid region are the compressible Navier-Stokes equation. In Cartesian coordinates and using the Einstein convention, the mass balance is expressed by:
\begin{equation*}
\frac{\partial \rho}{\partial t}+\frac{\partial \rho u_j}{\partial x_j}=0
\end{equation*}
where $x_j$, with $j=1,2,3$, are the Cartesian coordinates, $\rho$ is the fluid density and $u_j$ are the components of the velocity vector $\textbf{u}$. The momentum equation reads 
\begin{equation*}
\frac{\partial \rho u_i}{\partial t} + \frac{\rho u_{i}u_{j}}{\partial x_j}=-\frac{\partial p_{rgh}}{\partial x_{i}}-\frac{\partial \rho g_j x_j}{\partial x_i}+\frac{\partial}{\partial x_j}\left( \tau_{ij} + \tau_{tij}\right)
\end{equation*}
where $g$ the gravitational acceleration, $p_{rgh}=p-\rho g_j x_j$ the non-hydrostatic pressure, and $\tau_{ij}$ and $\tau_{tij}$ are the viscous and turbulent stresses, respectively.
The energy equation, written in terms of the specific enthalpy $h$, is expressed as:
\begin{equation*}
    \frac{\partial \rho h}{\partial t}+\frac{\partial}{\partial x_j}\left(\rho u_j h\right)+\frac{\partial \rho k}{\partial t}+\frac{\partial}{\partial x_j}\left(\rho u_j k_e \right)=-\frac{\partial q_i}{\partial x_i}+\frac{\partial p}{\partial t}-\rho g_j u_j +\frac{\partial}{\partial x_j}\left(\tau_{ij}u_i\right) + \dot{q}_s
\end{equation*}
in which $h=e+p/\rho$, where $e$ denotes the specific internal energy. Moreover we denote the kinetic energy $k_e=(u_i u_i)/2$, and $q_i$ the total heat transferred to the fluid, whereas $\dot{q}_s$ denote the heat source, including radiation (see \cite{bergman2011} for additional details).

For the solid region, only the energy equation is considered. The equation relates the rate of change of the solid enthalpy to the divergence of the heat conducted through the solid, and in Cartesian coordinate is given by:
\begin{equation}
\frac{\partial \rho h}{\partial t}=\frac{\partial}{\partial x_j}\left(\alpha \frac{\partial h}{\partial x_j}\right)    
\end{equation}
where  $\alpha=k/c_p$ denotes the thermal diffusivity, defined as the ratio between the thermal conductivity $k$ and the specific heat capacity $c_p$.
In order to guarantee a correct coupling between phases, at the interface between fluid and solid the temperature has to be the same, and the heat flux entering in one region at one side must be equal to the heat flux leaving the adjacent region, namely:
\begin{equation}
T_f=T_s, \quad \quad \quad
k_f\frac{dT_f}{dn}=-k_s\frac{dT_s}{dn},
\end{equation}
where $T_f$ and $T_s$ denote the fluid and solid temperatures, respectively, $n$ the interface normal, and $k_f$ and $k_s$ are the thermal conductivities for the fluid and solid regions, respectively.
\subsection*{Thermophysical properties}
For both fluid and solid regions the dependence of the thermophysical properties on temperature has been taken into account. The  heat  capacity $c_p$ for gas and air has been modeled according to the Janaf polymomial as $c_p=R(a_4T^4+a_3T^3+a_2T^2+a_1T+a_0)$, where $T$ is the temperature, $R$ is specific gas constant (287.05 [J/kg K] for air and 299.25 [J/kg K] for gas) and the coefficients $a_n$ for the two fluids are $a_{0,a}$=3.568, $a_{1,a}$=6.787e-4, $a_{2,a}$=1.554e-6, $a_{3,a}$=-3.299e-12, $a_{4,a}$=-4.664e-13, and $a_{0,g}$=3.569, $a_{1,g}$=5.473e-4, $a_{2,g}$=7.858e-6, $a_{3,g}$=-5.660e-10, $a_{4,g}$=1.296e-13.

The dynamic viscosity for air and gas has  been  modeled following the 
Sutherland law $\mu = A_s \sqrt{T}/(1+Ts/T)$, where $A_s$ and $T_s$ denotes the Sutherland's coefficient, defined for air and gas as follows: $A_{s_{air}}=1.458 \times 10^{-6}$ kg/m sec K$^{1/2}$, $T_{s_{air}}=110.4$ K; $A_{s_{gas}}=1.5544 \times 10^{-6}$ kg/m sec K$^{1/2}$, $T_{s_{gas}}=223.9$ K.
Finally, a cubic polynomial dependence has been employed for the solid thermal conductivity $k_s=bT^3+cT^2+dT+e$,
with $b=-5.0819 \times 10^{ -23}$, $c=1.1454 \times 10^{-19}$, $d=0.0113$, $e=10.004$.

\section{Results}
\label{num-res}

In this section we present an overview of the results of the numerical simulation campaign that has been carried out to define the optimal design of the new checkerboard recuperator. Our goal is to identify the geometrical configuration which maximizes the heat exchange, while keeping acceptable pressure drop values. For the different channel configurations considered (4x24, 5x24, 5x30, and 6x24) we have first investigated the effect of the wavy deformation and then analysed the effect of the radial scaling. Finally, the results have been compared with those related to a standard finned recuperator by means of the $\epsilon$-NTU method. This activity has been carried out in the framework of the European Project \textit{Burner 4.0}, funded by the Research Fund
for Coal and Steel, and the final design will be referred to as \textit{Burner 4.0} recuperator.

\subsection {The computational setup}
\label{comp-set}

All the numerical simulations have been performed using the open-source finite volume library OpenFOAM \cite{Weller1998}.  The flow has been considered steady and laminar and the solution strategy follows the conjugate heat transfer approach by means of the \texttt{chtMultiRegionFoam} solver. The total flow rate for the air and gas channels has been imposed to $\dot{Q}_a$= 0.01048 m$^3$/sec and $\dot{Q}_g$= 0.0464 m$^3$/sec, respectively, corresponding to $\dot{Q}_a$= 34 Nm$^3$/h and $\dot{Q}_g$= 37 Nm$^3$/h. The air flow enters the domain with an inlet temperature of Tin$_{air}$= 30$^\circ$C, whereas the gas inlet temperature is Tin$_{gas}$= 960 $^\circ$C. The overall length of the recuperator, including inlet-outlet manifolds is $L=$ 550 mm, whereas the external and internal diameter are $\phi_{ext}=$ 149 mm and $\phi_{int}=$ 65 mm respectively. In order to reduce the computational effort, exploiting the angular periodicity of the device, the computational domain consists of a single azimuthal element for air and gas regions, as shown in Figure \ref{fig:geom}. The cylindrical surfaces corresponding to the external and internal diameter have been considered adiabatic, and cyclic boundary conditions have been imposed on the lateral boundaries, for air, gas and solid regions. All the numerical studies have been performed on a computational grid of about 3.0 million of elements. 

As performance indicators, we considered the mean outlet temperatures for air and gas defined as the mass flow weighted average on the outlet sections, namely:
 \begin{equation}
          T_{out} = \frac{\int_{A_{out}}\rho u T \, dA}{\int_{A_{out}} \rho u \, dA},
          \label{eq:tout}
 \end{equation}
and the mean pressure drop $\Delta p = \frac{\Delta p_{air} + \Delta p_{gas}}{2}$, that is the average of the air and gas pressure drops between inlet and outlet

These performance indicators
 have been summarized in Table \ref{tab:baseline}, for all the Baseline cases with different number of channels. 
As expected the best performing configurations are those having more channels, in which the exchange surface is larger. In particular, the Baseline configuration with 6x24 channels is the one which guarantees the highest outlet temperature (689.3 $^{\circ}$C), still keeping the pressure drop to an acceptable value (137.4 Pa). In the following sections, we will investigate at which extent the wavy deformation (by increasing the exchange surface with respect to the baseline configuration) and the scaling deformations (by allowing a more uniform distribution of the flowrate through the channels) may improve the overall performance of the recuperator. 

\begin{table}[t]
\centering
\begin{tabular}{| l | c | c | c | c | }
 \hline
 & {4x24} & {5x24} & {5x30} &  {6x24} \\
\hline
{T$_{out AIR}$ [$^{\circ}$C]} & 590.3 & 633.7 & 654.4  & 689.3\\ 
{T$_{out GAS}$ [$^{\circ}$C]} & 531.3 & 594.5 & 475.5 & 447.2 \\
{${\Delta p}$ [Pa]}           & 73.7  & 99.1  & 111.9 & 137.4\\
\hline
\end{tabular}
 \caption{Performance analysis for the Baseline cases.}
 \label{tab:baseline}
\end{table}

\subsection{Wavy transformation}
\begin{figure}[h!]
    \centering
    \includegraphics[width=\textwidth]{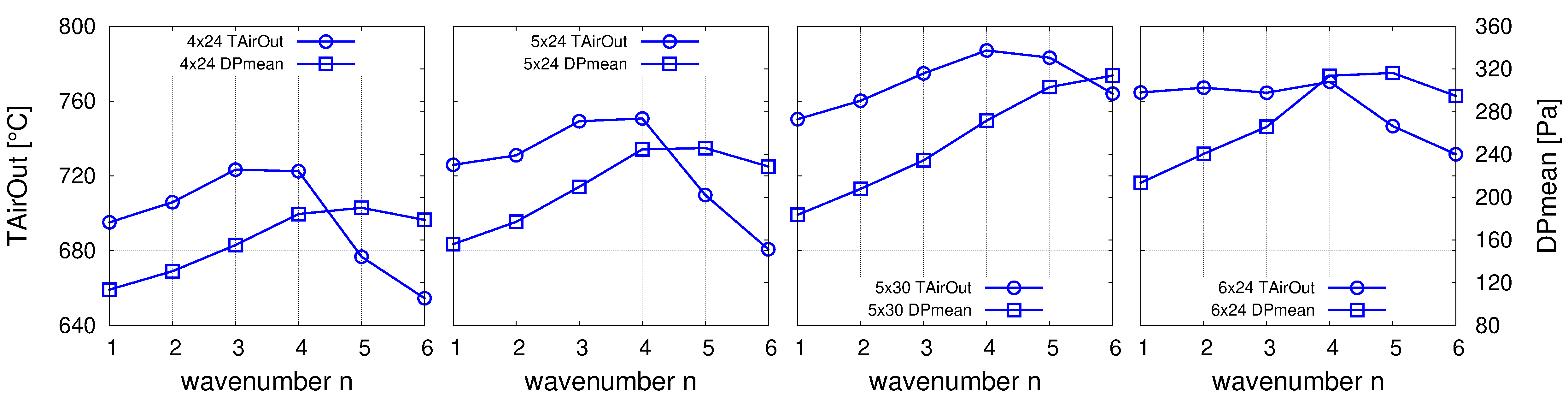}
    \caption{Behavior of the outlet air temperature and the mean pressure drop vs the wavenumber $n$, for the different channel configurations.}
    \label{fig:wavy}
\end{figure}

We first consider the wavy deformation defined in \eqref{eq:wavy} with the wavenumber $n$ ranging from 1 to 6 (integer wavenumbers have been considered to ease the assembling between the core sections and the end and middle collectors).  The results are collected in Table \ref{tab:wavy} for the mean air outlet temperature and mean pressure drop for different values of wavenumber, with the better performing configurations highlighted in bold font. 
\begin{table}[b]
\centering
\begin{tabular}{ | c | l | c | c | c | c | c | c | }
 \hline
{Channels} & {Wavenumber n} & {1} & {2} & {3} & {4} & {5} & {6}  \\
\hline
 & {T$_{out AIR}$ [$^{\circ}$C]} & 695.1 & 705.9 & \textbf{723.3} & 722.5 & 676.7 & 654.6 \\ 
 4x24 & {T$_{out GAS}$ [$^{\circ}$C]} & 447.7 & 435.3 & {405.3} & 414.9 & 453.2 &  479.5\\
 & {${\Delta p}$ [Pa]} & 113.5 & 130.8 & 155.2 & 184.3 & 190.1 & 178.8 \\
 \hline
 & {T$_{out AIR}$ [$^{\circ}$C]} & 725.9 & 731.0 & 749.2 & \textbf{750.7} & 709.7 & 680.7\\ 
 5x24 & {T$_{out GAS}$ [$^{\circ}$C]} & 417.9 & 409.4 & 396.0 & {391.0} & 422.7 & 454.7\\
 & {${\Delta p}$ [Pa]} & 156.1 & 177.0 & 209.8 & 244.8 & 246.0 & 228.7 \\
 \hline
 & {T$_{out AIR}$ [$^{\circ}$C]} & 750.3 & 760.2 & 774.8 & \textbf{787.0} & 783.2 & 764.0\\ 
 5x30 & {T$_{out GAS}$ [$^{\circ}$C]} & 394.3 & 385.4 & 373.3 & 350.2 & 361.3 & 382.1 \\
 & {${\Delta p}$ [Pa]} & 183.5 & 207.8 & 234.5 & 271.7 & 303.1 & 313.9 \\
\hline
 & {T$_{out AIR}$ [$^{\circ}$C]} & 764.5 & 767.1 & 764.4 & \textbf{770.3} & 746.6 & 731.5\\ 
 6x24 & {T$_{out GAS}$ [$^{\circ}$C]} & 380.3 & 380.1 & 376.4 & {372.0} & 389.9 & 411.5\\
 & {${\Delta p}$ [Pa]} & 213.5 & 240.6 & 265.9 & 313.7 & 316.3 & 294.8\\
 \hline
\end{tabular}
 \caption{Parametric analysis for the wavy deformation (in bold the best performance for each channel configuration)}
 \label{tab:wavy}
\end{table}
Note that the maximum value of the mean outlet air temperature for the 4x24 configuration, corresponding to 723.3 $^\circ$C, occurs at wavenumber $n=3$. For the other channel configurations, the best performances are achieved with $n=4$. The highest mean outlet air temperature (787.0 $^\circ$C) is obtained for the 5x30 configuration, with a 20\% improvement with respect to the corresponding baseline configuration. Although, as expected, the pressure jump increases when the wavy transformation is applied, the values obtained for all the configurations are always acceptable. 
The same results are also displayed in Figure \ref{fig:wavy}, where the trends of the mean outlet air temperature and mean pressure drop w.r.t. the wavenumber.

\vspace{-2mm}
\subsection{Scaling transformation}
Starting from the best performing configurations obtained with the wavy transformation, a second parametric analysis has been performed varying the scaling parameter $q$. Different ranges on $q$ have been considered to guarantee for each channel configuration a minimal radial thickness. The results are collected in Table \ref{tab:scaling}, with the better performing configurations highlighted in bold font. 
\begin{table}[b]
\centering
\begin{tabular}{| c | l | c | c | c | c | c | c | c | c | c |}
\hline
& {Scaling $q$} & {0.7} & {0.75} & {0.8} & {0.85} & {0.9} & {0.95} & {1.0} & {1.05} & {1.1}  \\
\hline
 & {T$_{out AIR}$ [$^{\circ}$C]} & 702.3 & 714.2 & 721.8 & 720.9 & \textbf{726.6} & 723.3 & 722.6 & 721.3 & 716.9\\
 4x24 ($n=3$)  & {T$_{out GAS}$ [$^{\circ}$C]} & 426.7 & 422.9 & 428.3 & 420.0 & 413.9 & 407.7 & 405.3  & 402.7 & 412.7\\
 & {${\Delta p}$ [Pa]} & 164.8 & 163.9 & 162.6 & 161.2 & 160.8 & 156.9 & 155.2 & 151.2 & 150.5\\
 \hline
 & {T$_{out AIR}$ [$^{\circ}$C]} & & & 736.8 & 748.8 & 749.8 &  \textbf{755.3} & 750.7 & 742.4 &  722.3 \\
5x24 ($n=4$) & {T$_{out GAS}$ [$^{\circ}$C]} & & & 403.0 & 391.2 & 386.3 & 388.4 & 391.0 & 403.4 & 421.5\\
 & {${\Delta p}$ [Pa]} & & & 246.0 & 245.4 & 244.9 & 246.4 & 244.8 & 242.0 & 239.7 \\
 \hline
 & {T$_{out AIR}$ [$^{\circ}$C]} & & & 776.7 & 785.0 & 786.6 &  \textbf{790.9} & 787.0 & 783.6 &  776.0\\
 5x30  ($n=4$) & {T$_{out GAS}$ [$^{\circ}$C]} & & & 367.6 & 371.2 & 362.1 & 358.8 & 350.2 & 357.3 & 358.1 \\
 & {${\Delta p}$ [Pa]} & & & 284.6 & 280.7 & 278.5 & 276.5 &  271.7 & 267.5 & 259.1\\
 \hline
 & {T$_{out AIR}$ [$^{\circ}$C]} & & & & & 773.7 & \textbf{777.6} & 770.3 & 766.7 & 747.0 \\
 6x24  ($n=4$) & {T$_{out GAS}$ [$^{\circ}$C]} & & & & & 421.5 & 378.5 & 372.0 & 387.2 & 393.1 \\
 & {${\Delta p}$ [Pa]} & & & & & 318.8 & 316.3 & 313.7 & 307.4 & 298.5 \\
 \hline
\end{tabular}
 \caption{Parametric analysis for the scaling deformation (in bold the best performance for each channel configuration)}
 \label{tab:scaling}
\end{table}

Also in this case, the best performance has been achieved by the 5x30 configuration with a scaling factor $q=0.95$ increasing the mean air outlet temperature to 790.0 $^\circ$ C, corresponding to an increase of 21\% with respect to the baseline configuration, thus further improving, although to a limited extent, the overall performance.

\begin{figure}[h!]
    \centering
     \includegraphics[width=\textwidth]{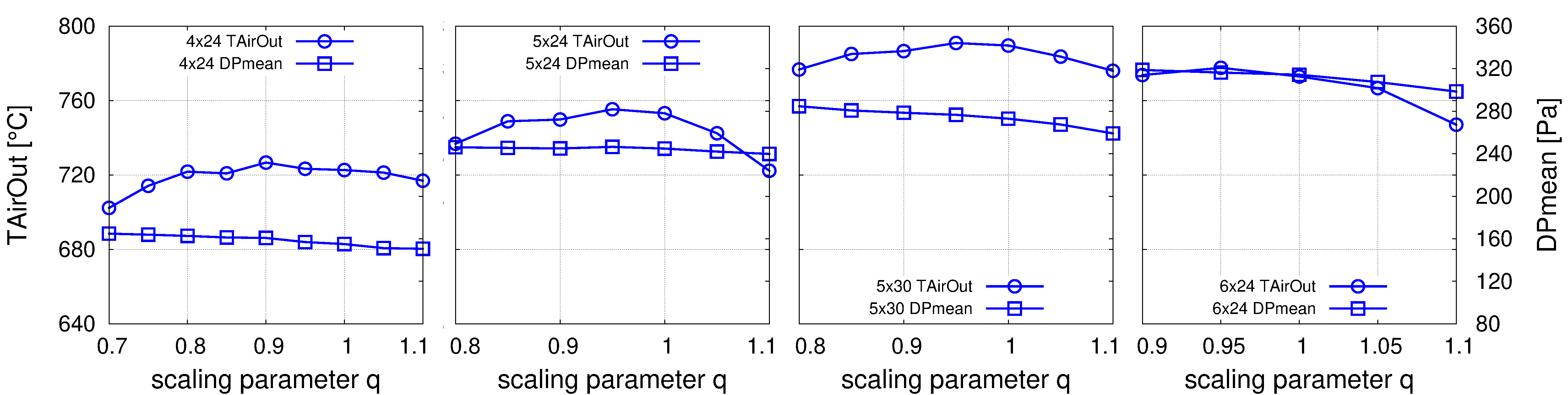}
    \caption{Behavior of the outlet air temperature and the mean pressure drop vs the scaling parameter $q$, for the different channel configurations.}
    \label{fig:perf_scal}
\end{figure}

In order to quantify the possible contribution of radiation, on this last configuration an additional simulation has been performed using the Finite Volume Discrete Ordinates Method (FvDOM) radiation model (see \cite{modest2013}  for more details). From the results reported in Table \ref{tab:radiation}, we observe that thermal radiation slightly affects the performance, allowing the outlet mean air temperature to increase of a few degrees Celsius. We point out that similar analyses (not reported here for brevity) have been carried out for the best configuration of each checkerboard layout, showing the same limited effect of the radiative contribution.

\begin{table}
\scriptsize
\centering
\begin{tabular}{| l | c | c | }
\hline
 & {no Rad} & {FvDOM Rad model}\\
\hline
{T$_{out AIR}$ [$^{\circ}$C]} & 790.9 & 793.2 \\ 
{T$_{out GAS}$ [$^{\circ}$C]} & 358.8 & 357.0 \\
{${\Delta p}$ [Pa]} & 276.5 & 276.6\\
\hline
\end{tabular}
 \caption{Comparison for the checkerboard layout 5x30, n=4, q=0.95 between the performance indicators with and without fvDOM radiation model.}
 \label{tab:radiation}
\end{table}

\vspace{-5mm}
\subsection{Finned recuperator}
\label{sec:finned}

\begin{figure}[h!]
    \centering
    \includegraphics[width=0.5\textwidth]{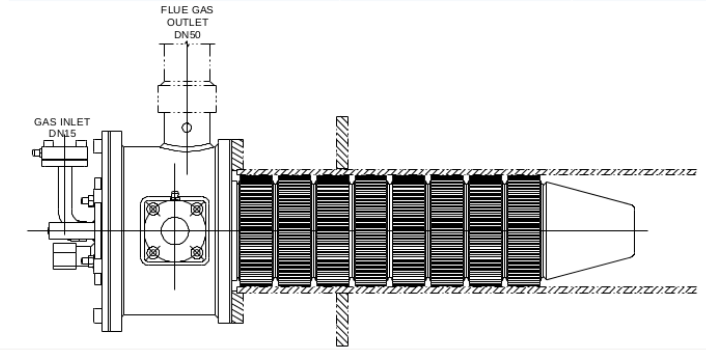}
    \hspace{5mm}
     \includegraphics[width=0.35\textwidth]{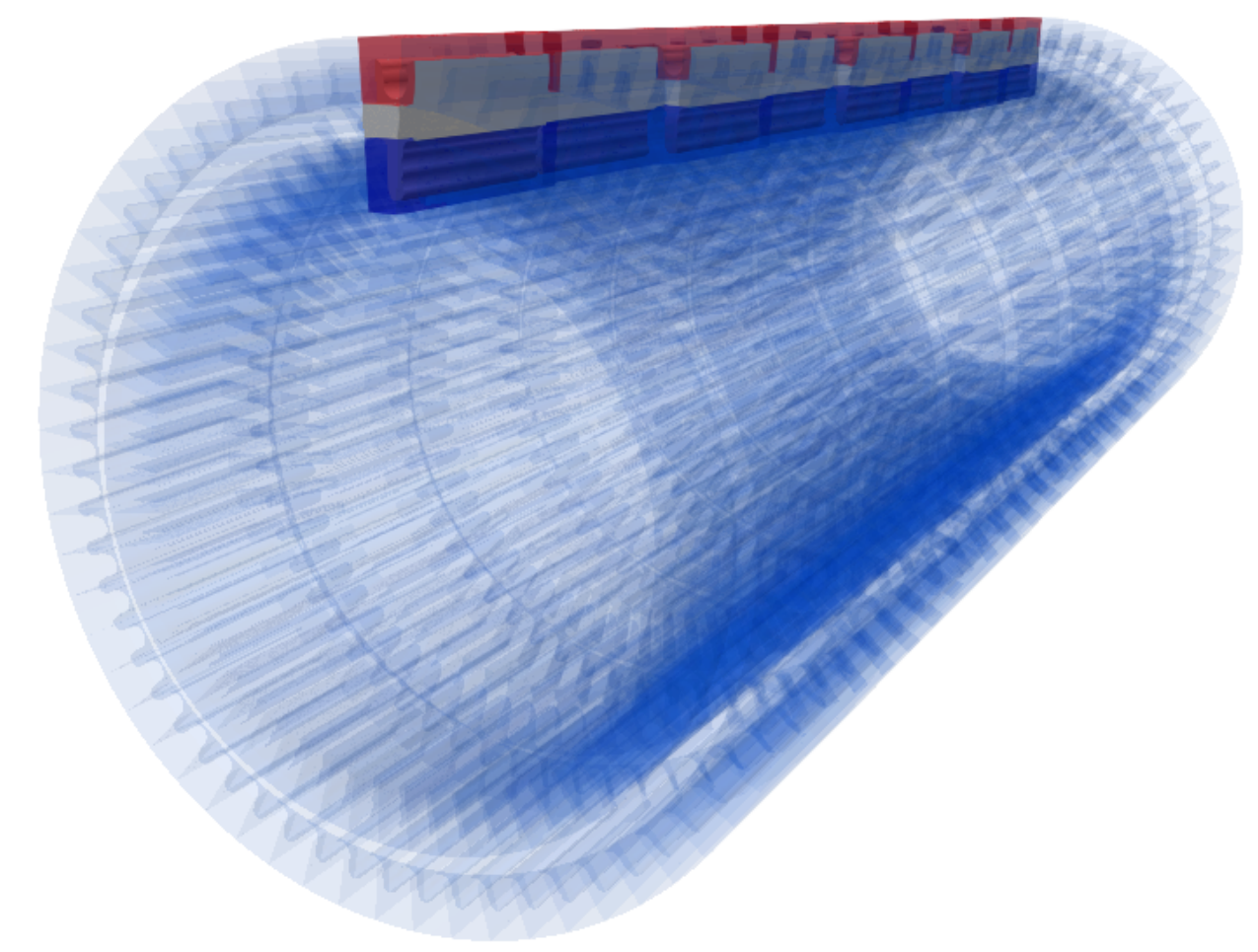}
     \caption{Standard finned recuperator geometry (left) and corresponding computational domain (right). Blue and red areas denote the air and gas regions, respectively.}
         \label{fig:finned}
\end{figure}

In order to quantify the potential benefit of the new checkerboard design, we have also considered, as a reference configuration, the standard finned recuperator, representing the existent technology currently used by Tenova, which is composed by 8 finned staggered element, with 60 radial fins, as shown in Figure \ref{fig:finned} (left).
The axial length is L= 400 mm, whereas the external and internal diameter are $\phi_{ext}$= 149 mm and $\phi_{int}$= 99 mm, respectively. The computational model has been developed considering a computational domain consisting of a single finned radial element, see Figure \ref{fig:finned} (right). 
The numerical tests have been performed, using the OpenFOAM \texttt{chtMultiRegionFoam} solver under the same temperature and flow rate condition reported in Section \ref{comp-set} for the checkerboard cases. The surfaces corresponding to the external and internal diameters have been assumed to be adiabatic, and cyclic boundary condition has been imposed on lateral walls.
The numerical simulations have been carried out on a computational grid composed by 9 million elements, with and without $k-\omega$ SST turbulence model. The effect of radiation has been also investigated by means of the fvDOM (Finite Volume Discrete Ordinates Method) radiation model. 
From the results reported in Table \ref{tab:fin}, in terms of mean outlet temperatures and mean pressure drop, we can note that the adoption of the $k-\omega$ SST turbulence model has a limited impact on the results, whereas, for this kind of configuration, it is fundamental to account for the radiative contribution.

\begin{table}
\centering
\begin{tabular}{| l | c | c | c | c | c | c |}
\hline
 &{Laminar} & {k-$\omega$ SST} & {Laminar + fvDOM} \\
\hline
{T$_{out AIR}$ [$^{\circ}$C]} & 569.2 & 563.7 & 637.8 \\
{T$_{out GAS}$ [$^{\circ}$C]} & 567.2 & 566.8 & 501.6\\
{${\Delta p}$ [Pa]} & 340.8 & 336.1 & 351.1 \\
 \hline
\end{tabular}
 \caption{Performance indicators for the finned recuperator simulations.}
 \label{tab:fin}
\end{table}

\vspace{-5mm}
\subsection{Performance analysis through the $\epsilon$-NTU method}
\label{perf-an}

\begin{figure}[b!]
    \centering
     \includegraphics[width=\textwidth]{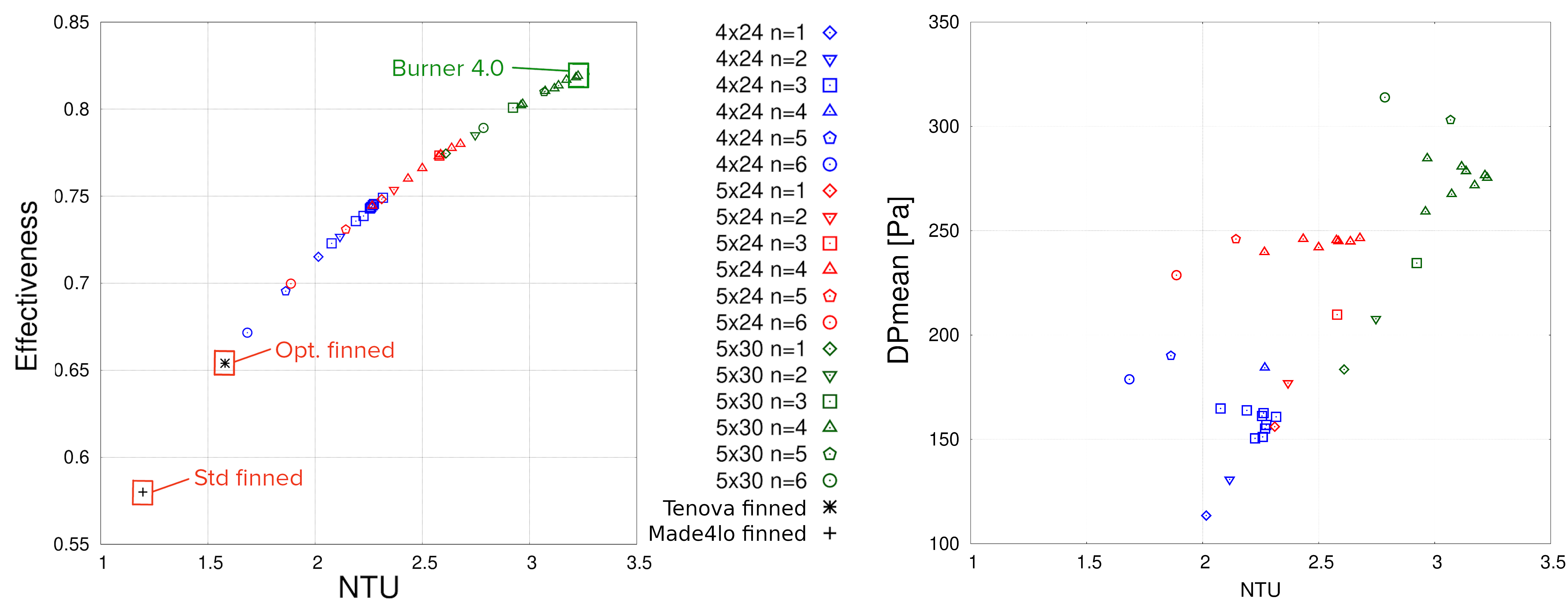}
     \caption{Effectiveness $\epsilon$ (left) and mean pressure drop DP (right) as a function of NTU for all the heat exchanger configurations \cite{made4lo}}
     \label{fig:ntu}
\end{figure}

In order to evaluate the thermal performance of the different analyzed configuration, a comparison with the $\epsilon$-NTU method \cite{shah2003} has been carried out. 
According to the $\epsilon$-NTU method, the effectiveness $\epsilon=\frac{q}{q_{max}}$ is the ratio between the actual heat transfer rate $q$ from the hot to the cold fluid and the maximum possible heat transfer rate $q_{max}$ that can be estimated based on general thermodynamical considerations. Considering a counterflow heat exchanger of infinite surface area, an overall energy balance for the two fluids can be written as
\begin{equation}
    q=C_h\left(T_{h,i}-T_{h,o}\right)=C_c\left(T_{c,o}-T_{c,i}\right)
\end{equation}
where $T_{h,i}$ and $T_{h,o}$ denotes the inlet and outlet temperature of the hot flow, respectively, whereas $T_{c,i}$ and $T_{c,o}$ stand for the inlet and outlet temperature of the cold flow. $C_h$ and $C_c$ are the heat capacity rate for the hot and cold fluid, defined as $C_h=\left(c_p \dot{m}\right)_h$ and $C_c=\left(c_p \dot{m}\right)_c$, where $c_p$ is the specific heat capacity rate and $\dot{m}$ the mass flow rate.
Over the infinite flow length, the hot fluid temperature will approach the inlet temperature of the cold fluid, resulting in $T_{h,o}$=$T_{c,i}$. The maximum possible heat transfer can be defined as
\begin{equation}
    q_{max}=C_{min}\left(T_{h,i}-T_{c,i} \right)=C_{min}\Delta T_{max}
\end{equation}
in which $C_{min}$=$C_c$ if $C_c$ $<$ $C_h$, and $C_{min}$=$C_h$ if $C_f$ $<$ $C_c$, and $\Delta T_{max}$ denotes the maximum temperature difference. It can be proven, see e.g. \cite{shah2003}, that $\epsilon$ can also be written as:
\begin{equation}
    \epsilon = \frac{UA}{C_{min}}\frac{\Delta T_m}{\Delta T_{max}}
\end{equation}
where $U$ is the overall heat transfer coefficient, and $A$ the heat transfer surface area. More generally, for each heat exchanger layout, the effectiveness $\epsilon = f\left(NTU, C^*\right)$ is a function of the number of transfer units $NTU=UA/C_{min}$, and the heat capacity rate ratio $C^*=C_{min}/C_{max}$. The NTU is a design parameter, denoted also as the non-dimensional heat transfer size of the exchanger. NTU provides a compound measure of the heat exchanger size through the product of the total heat transfer area $A$, and the overall heat transfer coefficient $U$. At low values of NTU, $\epsilon$ is low, whereas with increasing value of NTU $\epsilon$ increases, approaching a thermodynamic asymptotic values \cite{shah2003}.
Figure \ref{fig:ntu} (left) shows the values of $\epsilon$ and NTU for the CFD results of both checkerboard and finned layout. The label \textit{Tenova finned} denotes the standard finned recuperator analyzed in the previous section, whereas the label \textit{Made4Lo finned} refers to the CFD data provided by Tenova, w.r.t another type of finned recuperator developed in additive manufacturing in a previous research project (Made4Lo  \cite{made4lo}).
We observe that the points distribution in the (NTU,$\epsilon$) plane follows the asymptotic trend typical of counterflow heat exchanger layouts at constant $C^*$. 

Note that the lowest values of $\epsilon$ and NTU occur for both the finned exchanger, and are approximately 0.58 and 0.65 for $\epsilon$ and 1.2 and 1.6 for NTU, for the standard and Made4Lo finned recuperators, respectively. Compared to the finned cases, the checkerboard design shows higher values of $\epsilon$ and NTU, for each $n \times m$ configuration. According to the previous parametric analysis, the 5x30 checkerboard layout at $n$=4 and $q$=0.97 turns out to be the most promising configuration, with $\epsilon$ up to $\approx$ 0.82, and NTU up to $\approx$ 3.2. The behavior of $\epsilon$ and NTU in Figure \ref{fig:ntu} (left) highlights the effect of the total heat transfer area on the heat exchanger performance, and clearly suggests that the development of the checkerboard design based on the Tenova patent \cite{dellarocca2018b}, could be potentially very advantageous in order to improve the existing technology for recuperative burners. To complete the analysis, in Figure \ref{fig:ntu} (right) the behavior of the mean pressure drop DP versus NTU is given for each checkerboard configuration. We observe that the mean pressure drop is always within a range of admissibility. Note that, if the goal is NTU $>$ 3, the corresponding cost in terms of mean pressure drop is $>$ 250 Pa, especially the optimal configuration corresponds to a pressure drop DP $\approx$ 275 Pa.

\section{Experimental validation}\label{validation}
\begin{figure}[b!]
    \centering
     \includegraphics[width=0.565\textwidth]{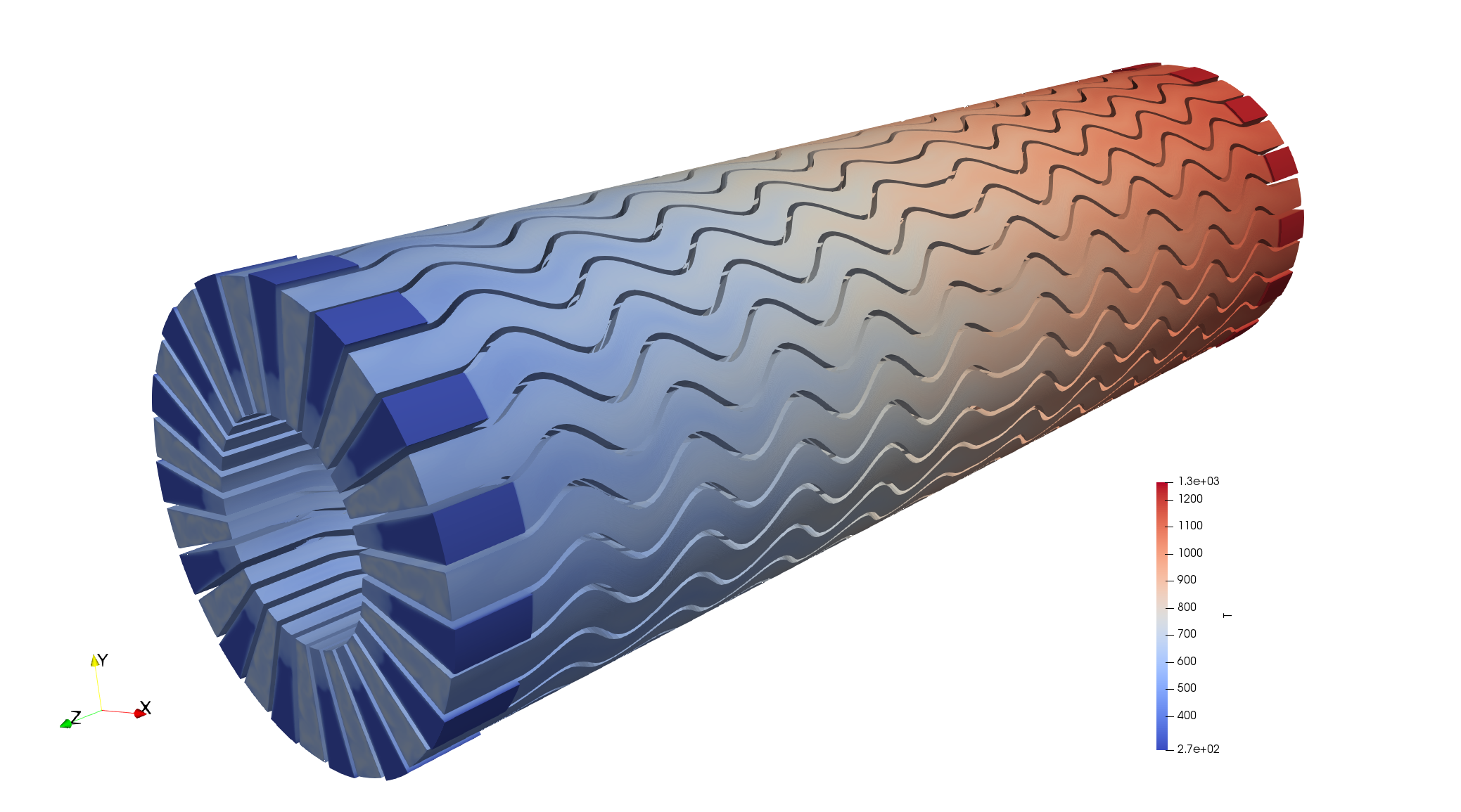}
     \includegraphics[width=0.425\textwidth]{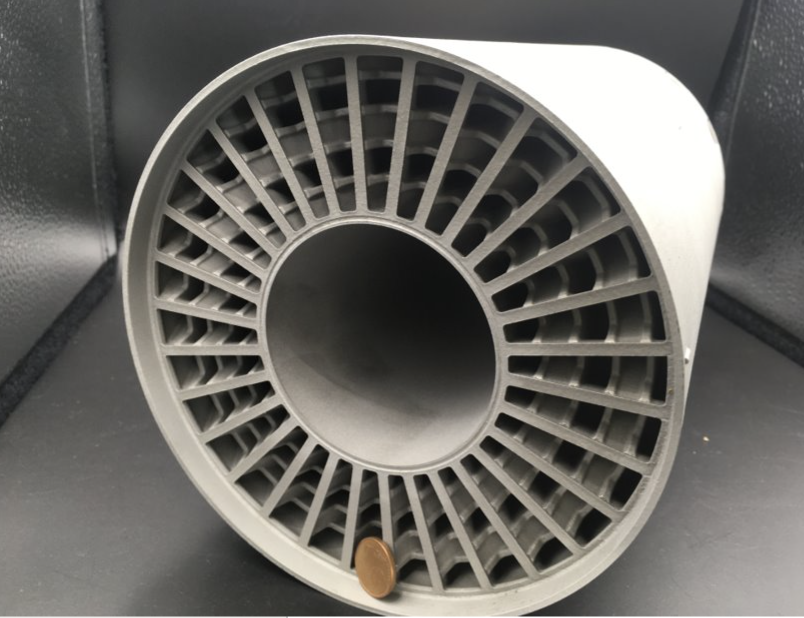}
     \caption{Final optimal checkerboard design (left) and 3D printed recuperator (right).}
     \label{fig:printed}
\end{figure}

The checkerboard recuperator with the optimal design as identified in the previous sections  (see Figure \ref{fig:printed}, left) has been 3D printed and assembled at the Department of Mechanics of Politecnico di Milano and experimentally tested by Rina S.p.a, also partner of the Burner 4.0 European project (see Figure \ref{fig:printed}, right). 

The results of the experimental tests carried out at RINA-CSM laboratory in Dalmine (Italy) comparing different recuperators at different inflow gas temperatures are reported in Table \ref{tab:exper}, in terms of outflow gas temperatures. The lower the outflow gas temperature value, the higher the recuperator efficiency.

The experimental campaign carried out by RINA-CSM had the main goal to verify the thermal efficiency of the Tenova self-recuperative burner equipped with different heat exchanger types. The self-recuperative prototype version selected for the experimental characterization were
\begin{itemize} 
\item the self-recuperative burner currently installed at a reference industrial plant (labeled \emph{Original});
\item the self-recuperative burner with standard heat exchangers
(labeled \emph{Tenova finned});
\item the self-recuperative burner with optimized heat exchanger (labeled \emph{Made4Lo finned});
\item the self-recuperative burner with checkerboard heat exchanger by additive manufacturing (labeled \emph{Burner4.0 checkerboard}).
\end{itemize}
The furnace used for the test campaign is composed of four identical modules. This design allows to reduce the furnace length in relation to the radiant tube thermal power. For the experimental test the 3+1 modules configuration was adopted. The modules are assembled by flanges with nuts and bolts. The modules are supported by a movable frame on wheels. The overall dimension of the furnace is 3.10 m in length and about 1.0 m x 1.0 m in square section. The internal section of the chamber is 0.5 m x 0.5 m, while the internal useful length with 3+1 modules is 2.08 m. The combustion system and the radiant tube were installed on the front wall of the furnace. The radiant tube used for the experimentation was of single end type and it was equipped with seven thermocouples on outer surface to evaluate the temperature distribution. Other thermocouples were used to measure the inlet and outlet temperature of flue gases through the recuperator, as well as the temperature of preheated combustion air exiting from the recuperator.

\begin{table}
\centering
\begin{tabular}{| c | c | c | c |}
\hline
{T$_\text{in, gas}$ [$^{\circ}$C]} & 850 & 950 & 1050 \\
\hline
Original & 597 & 620 & 647\\
Tenova finned & 443 & 452 & 497\\
Made4Lo finned & 462 & 481 & 510\\
Burner4.0 checkerboard & 326 & 331 & 357 \\
 \hline
\end{tabular}
 \caption{Mean gas outflow temperatures [$^\circ$C] for different recuperators at different gas inflow temperatures.}
 \label{tab:exper}
\end{table}

The superior performances of the optimized checkerboard configuration is clearly demonstrated by these results, which correspond to an efficiency $\varepsilon$ of the original reference recuperator of about $50\%$, an improved efficiency of the Tenova finned configurations between $60\%$ and $70\%$, while the checherboard configuration reaches an effienciency higher than $80\%$ as predicted by our CFD campaign. 

\section{Conclusion}
\label{conclusion}

In this paper we have presented an extensive numerical investigation carried out to optimize the design of the checkerboard air-gas recuperator based on the Tenova patent \cite{dellarocca2018b}. Different parametric studies have been performed by means of the CFD software OpenFOAM, for various geometrical configurations. The heat excharger efficiency has been assessed by means of the $\epsilon$-NTU method. The perfomances have been compared with those of traditional double pipe finned recuperator and significant enhancement of the exchanger efficiency has been observed adopting the new checkerboard design. 
The outcomes of our analysis suggest that the development of the present checkerboard heat exchanger could be potentially very advantageous in order to improve the existing technologies for heat recovery in high temperature industrial systems, as in burners for industrial furnaces.
\section*{Acknowledgments}

This project has received funding from the Research Fund for Coal and Steel under grant agreement No 847237.
The authors wish to express their thanks to Barbara Previtali, Leonardo Caprio and Andrea Valensin (DMEC, Politecnico di Milano) and Guido Jochler and Irene Luzzo (RINA-CSM) for their invaluable contribution in the design definition and experimental campaign. The authors gratefully acknowledge the computing time obtained by the CINECA Supercomputing Center with the Iscra-C StoheCFD and SoheCFD2 projects. 

\centering
\vspace{3mm}
 \includegraphics[width=0.3\textwidth]{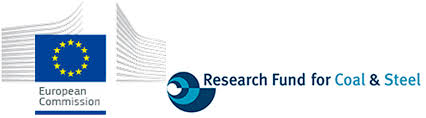}


\vspace{-3mm}

\bibliographystyle{ieeetr}
\bibliography{bibliography}

\end{document}